\pdfoutput=1
\RequirePackage{ifpdf}
\ifpdf 
\documentclass[pdftex]{sigma}
\else
\documentclass{sigma}
\fi

\begin{document}

\allowdisplaybreaks

\renewcommand{\thefootnote}{}

\renewcommand{\PaperNumber}{023}

\FirstPageHeading

\ShortArticleName{On the 80th Birthday of Dmitry Borisovich Fuchs}

\ArticleName{On the 80th Birthday of Dmitry Borisovich Fuchs\footnote{This paper is a~contribution to the Special Issue on Algebra, Topology, and Dynamics in Interaction in honor of Dmitry Fuchs. The full collection is available at \href{https://www.emis.de/journals/SIGMA/Fuchs.html}{https://www.emis.de/journals/SIGMA/Fuchs.html}}}

\Author{Alice FIALOWSKI~$^{\dag^1}$,
Ekaterina FUCHS~$^{\dag^2}$,
Elena FUCHS~$^{\dag^3}$,
Boris KHESIN~$^{\dag^4}$,\newline
Alexandre KIRILLOV~$^{\dag^5}$,
Fedor MALIKOV~$^{\dag^6}$,
Valentin OVSIENKO~$^{\dag^7}$,\newline
Alexei SOSSINSKY~$^{\dag^8}$
and Serge TABACHNIKOV~$^{\dag^9}$}

\AuthorNameForHeading{B.~Khesin, F.~Malikov, V.~Ovsienko, S.~Tabachnikov}

\Address{$^{\dag^1}$~Institute of Mathematics, University of P\'ecs, P\'ecs, Hungary}
\EmailDD{\href{mailto:fialowsk@ttk.pte.hu}{fialowsk@ttk.pte.hu}}
\URLaddressDD{\url{http://web.cs.elte.hu/~fialowsk/}}

\Address{$^{\dag^2}$~Department of Mathematics, City College of San Francisco, CA 94112, USA}
\EmailDD{\href{mailto:efuchs@ccsf.edu}{efuchs@ccsf.edu}}
\URLaddressDD{\url{https://www.ccsf.edu/Info/Faculty_In_Review/7272/}}

\Address{$^{\dag^3}$~Department of Mathematics, UC Davis, One Shields Ave Davis, CA 95616, USA}
\EmailDD{\href{mailto:efuchs@math.ucdavis.edu}{efuchs@math.ucdavis.edu}}
\URLaddressDD{\url{https://www.math.ucdavis.edu/~efuchs/}}

\Address{$^{\dag^4}$~Department of Mathematics,
University of Toronto, Toronto, ON M5S 2E4, Canada}
\EmailDD{\href{mailto:khesin@math.toronto.edu}{khesin@math.toronto.edu}}
\URLaddressDD{\url{http://www.math.toronto.edu/khesin/}}

\Address{$^{\dag^5}$~Department of Mathematics, University of Pennsylvania, Philadelphia, PA 19104-6395, USA}
\EmailDD{\href{mailto:kirillov@math.upenn.edu}{kirillov@math.upenn.edu}}
\URLaddressDD{\url{https://www.math.upenn.edu/~kirillov/}}

\Address{$^{\dag^6}$~Department of Mathematics, University of Southern California,
Los Angeles, CA 90089, USA}
\EmailD{\href{mailto:fmalikov@usc.edu}{fmalikov@usc.edu}}
\URLaddressD{\url{https://dornsife.usc.edu/cf/faculty-and-staff/faculty.cfm?pid=1003489}}

\Address{$^{\dag^7}$~CNRS, Laboratoire de Math\'ematiques de Reims,
51687 Reims cedex 2,
France}
\EmailDD{\href{mailto:valentin.ovsienko@univ-reims.fr}{valentin.ovsienko@univ-reims.fr}}
\URLaddressDD{\url{http://ovsienko.perso.math.cnrs.fr}}

\Address{$^{\dag^8}$~Independent University of Moscow,
Bolshoy Vlasyevskiy Pereulok 11, 119002, Moscow, Russia}
\EmailDD{\href{mailto:asossinsky@yandex.ru}{asossinsky@yandex.ru}}
\URLaddressDD{\url{https://users.mccme.ru/abs/}}

\Address{$^{\dag^9}$~Department of Mathematics,
Pennsylvania State University, University Park, PA 16802, USA}
\EmailDD{\href{mailto:tabachni@math.psu.edu}{tabachni@math.psu.edu}}
\URLaddressDD{\url{http://www.personal.psu.edu/sot2/}}

\ArticleDates{Received February 29, 2020; Published online April 07, 2020}

\Abstract{This article is a collection of several memories for a
special issue of SIGMA devoted to Dmitry Borisovich Fuchs.}


\Classification{01A60; 01A61}

\newpage

\renewcommand{\thefootnote}{\arabic{footnote}}
\setcounter{footnote}{0}

\section[Preface by B.~Khesin, F.~Malikov, V.~Ovsienko and S.~Tabachnikov]{Preface by B.~Khesin, F.~Malikov, V.~Ovsienko\\ and S.~Tabachnikov}

This special issue of SIGMA is devoted to Dmitry Borisovich Fuchs, an outstanding mathematician and our teacher, on the occasion of his 80th birthday. It comprises many papers of his friends and colleagues,
reflecting broad mathematical interests of D.B.\ himself.

Below we collected several memories: those of his daughters Katia and Lyalya, of his old-time friends
A.~Sossinsky and A.~Kirillov, and his students and colleagues A.~Fialowski and V.~Ovsienko.

There are many (semi-legendary) stories about Fuchs and, by way of introduction, we present a few specimen of this lore.

\begin{figure}[ht]\centering
\includegraphics[width=4.5in]{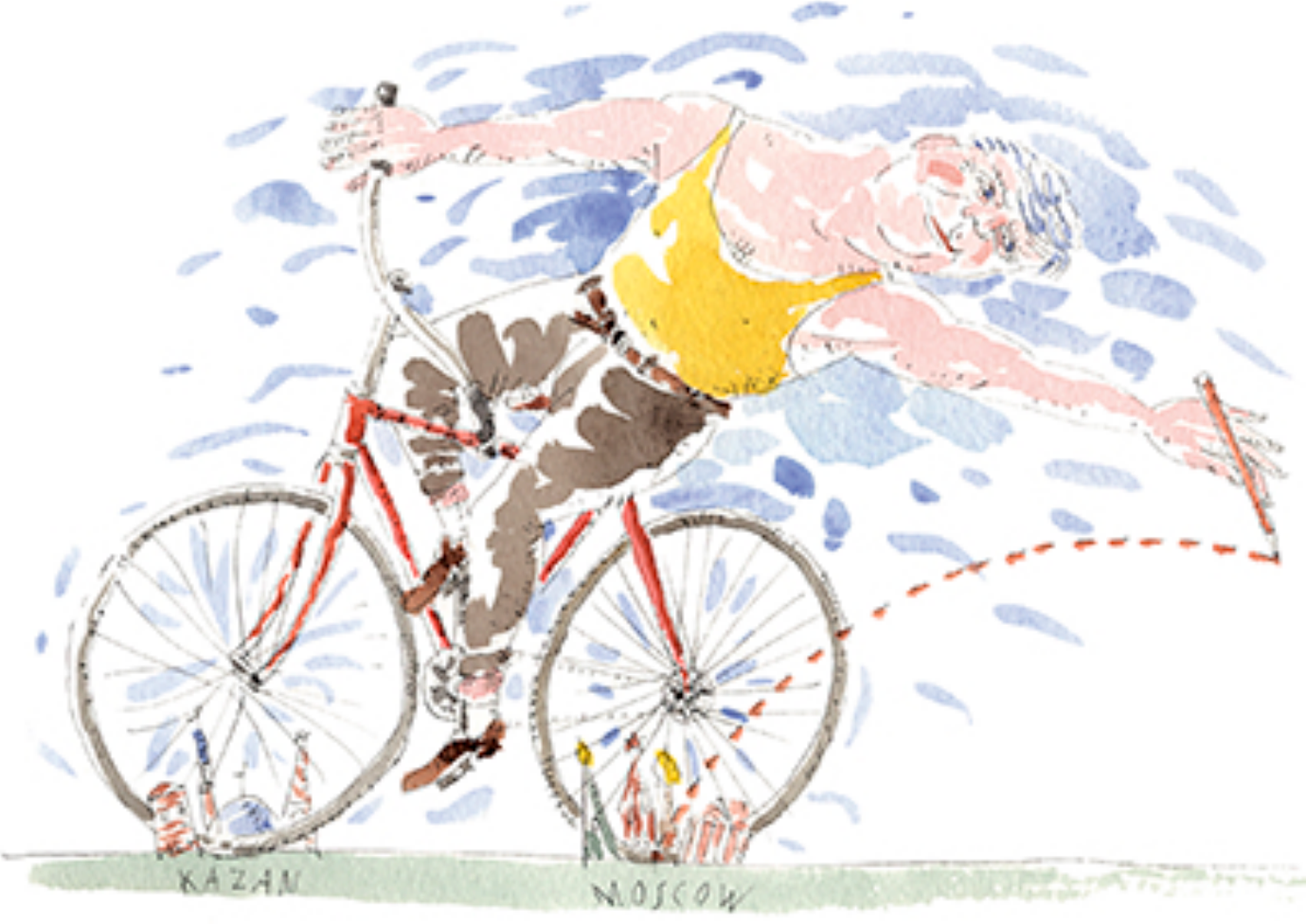}
\caption{Fuchs is an avid bicyclist, making very long distance trips, in particular, from Moscow to Kazan (about 500 mi). This picture was drawn by Sergei Ivanov, the artist who illustrated the book ``Mathematical Omnibus'' by D.~Fuchs and S.~Tabachnikov, AMS, Providence, RI, 2007.}
\end{figure}

In his memories of Arnold's seminar in Moscow D.B.~Fuchs writes (see~\cite{Arn})
``My role there was well
established: I had to resolve any topology-related difficulty. Some of my friends
said that at Arnold's seminar I was a `cold topologist'. Certainly, a non-Russian speaker
cannot understand this, so let me explain. In many Russian cities there
were `cold shoemakers' in the streets who could provide an urgent repair to your
footwear. They sat in their booths, usually with no heating (this is why they
were `cold'), and shouted, `Heels! \dots\ Soles! \dots' So I appeared as if sitting in a
cold booth and yelling, `Cohomology rings! \dots\ Homotopy groups! \dots\ Characteristic
classes! \dots'\,''

The topological reputation of Fuchs extended far beyond Arnold's seminar.
There was a~famous story where a group of mathematicians, ``all the best
people'' as Fuchs says, helped one of them, A.A.~Kirillov, to
move to his new apartment on the 14th floor of a 16-story building. It turned out that a
huge custom-made bookcase, which would occupy the entire wall of a
room, did not fit into the elevator. It had to be carried to the 14th floor
via the staircase but, even more importantly, it was so large that there
was a unique way of moving it through each flight of stairs, with the only
possible turn of the bookcase fitting the stairwell. (One of the participants was Ya.G.~Sinai,
who first formulated the problem as ``keep this damned bookcase connected'' and later ``keep the number of its components not greater than two''.)
 Eventually, when they got to the
apartment, there was a unique way to bring this bookcase into the
room, with no way to turn it around or change its orientation in any other way.

And once inside, it turned out
that the bookcase fit the room in the only possible position, with
the shelves facing the wall and the backside facing the room!
The only way to fix this was to carry the case all 14 flights down the
stairs and get it out to the street.
Then Fuchs, as the leading topologist, was given the task of finding the
correct ``initial conditions'' in order to start the process anew. He thought about it, gave the instructions
how to turn it, they did all the lifting via 14 flights of stairs up~-- and
the bookcase indeed fit in the room perfectly!

S.~Tabachnikov recalls (see~\cite{Tab}) ``At our first meeting DB asked me what I knew in mathe\-matics, and I started to name various books I had read (after I~mentioned Spanier's ``Algebraic topology'' DB remarked, characteristically, that it was the first time he met someone who had managed to finish the book). DB asked
a few questions (e.g., which simply connected 4-dimensional manifold I knew), and
my answers revealed the obvious: my extensive reading was almost a total waste
of time, I did not really understand the material. This was a devastating discovery
for me, and I mumbled something like: ``I should probably quit mathematics, I've
already lost too much time!'' This made DB laugh. ``You are 20, aren't you? It is
not too late to start, even if you had never read anything but `Murzilka'\footnote{A Russian magazine for preschoolers.} before''.

F.~Malikov remembers: Fuchs was giving a course on linear algebra and once, trying to entertain us, he devoted a lecture to a proof
of Abel's theorem (impossibility of solving equations of degree five in radicals). While proving it, he constructed a certain permutation. Then unexpectedly he turned to A.~Kanel', who was the indisputable genius of the class, and asked: ``Alesha, what could you say about this permutation?'' Kanel' replied: ``This permutation is even''. Fuchs said: ``Correct, but it's not its main property -- it is the trivial permutation''.

\begin{figure}[ht]\centering
\includegraphics[width=3in]{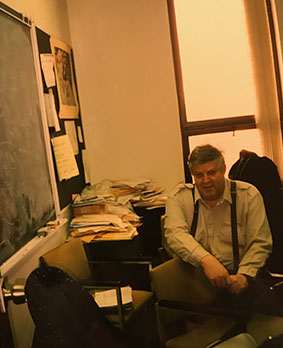}
\caption{In the office.}
\end{figure}

Yet another memory: Fuchs was giving a talk at USC. At some point, he said that, for a~sufficiently small epsilon, such and such held. One of the attendees, who was drowsing, suddenly woke up and asked how small this epsilon must be. Fuchs replied immediately: ``one tenth of an inch''.

B.~Khesin remembers the following story. Fuchs recalled that in the 70--80s he was at one of A.~Fomenko's public lectures on ``The New Chronology'', which was a ``hot'' subject at the time. At some point, when comparing time patterns related to two sequences of the Roman emperors, Fomenko said that the probability that such patterns were independent and not copied from one another was one out of a million. At this moment Fomenko tried to clarify his point: ``In order to see how negligible this probability of one millionth is, just imagine that you put a kettle on a hot stove, but instead of boiling the water freezes up! This is how rare such events are!'' Fuchs said that once he heard that explanation, he thought: ``Hmm \dots\ Moscow's population is 8 million \dots\ So every morning 8 people in this city put their kettles on the stove and have the water in the kettle frozen?!''

\section{Elena Fuchs}

\AuthorNameForHeading{El.~Fuchs}

\looseness=-1 One of my first memories of my father was waking up in the morning to find a sheet of math problems, handwritten by him and often involving me and my sister, in my room. My father would be gone to work by this point, but it was like a little piece of him that he left for me to enjoy.

\begin{figure}[ht]\centering
\includegraphics[width=3.5in]{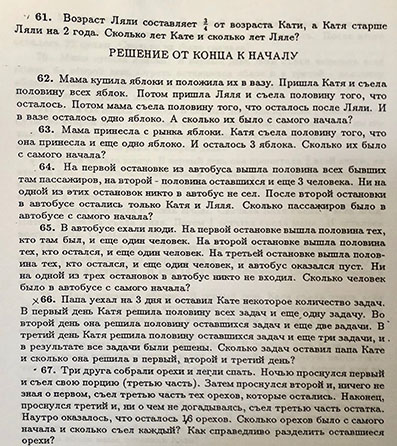}
\caption{A page from the problem book that Fuchs wrote for his daughters.}
\end{figure}

When I was a little older, after we moved to America, these sheets of problems came illustrated with a sketch of me doing something from one of the problems (I remember one picture in particular where I was diving into a swimming pool, something that he had taught me to do, in fact). Not everyone knows that my father is a wonderful artist: our apartment in Moscow had several of his oil paintings hanging on the walls. He doesn't paint much anymore, or talk about it, but I think he painted those paintings at a sad time in his life when the art helped him escape, so perhaps that is a good thing.

In fact, I do not remember a time in my life when math did not play a big role in our relationship. While we still lived in the USSR, I~remember that he was often away on what I~imagined to be an exciting adventure off in some other country. He would always come home with all sorts of treasures from the West that were unattainable at home at the time: colorful skirts, stuffed animals that looked like the real thing, fancy new pencils.

\begin{figure}[ht]\centering
\includegraphics[width=3in]{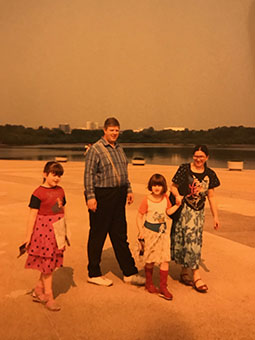}
\caption{The Fuchs family.}
\end{figure}

While my sister and I and the other kids in our neighborhood lounged in the sun on a mattress left on a ping-pong table at our dacha (summer home), my father would have frequent math visitors, and sometimes he would take me from there into town on his bike so that I could have ice cream while he met with yet another mathematician. Even when my sister and I did not fully appreciate it, our lives were engulfed in math. Wherever my father took us: to the skating rink, to a lesson of some kind, he always had his yellow pad of paper on which he would write and draw what back then looked like amusing doodles. During a month-long stay at St.~John's College in Cambridge in the mid-nineties, we would find his math toys at our apartment when my parents were away and play ridiculous games with them: one time we found a curvy ruler for tracing curves that one could shape into any curve one wanted and made it into an S-shape to advertise Safeway, the American grocery store, out of our apartment window. I~am sure the Brits were amused.

Rooms with floors covered in chalk dust have been places in which I have felt at home from a very young age, and it is all due to him. Indeed, the first time I~remember anyone asking me what I~wanted to be when I grew up, ``mathematician'' was always on the long list of professions I planned to have (along with opera singer, ballerina, and flautist \dots\ those did not pan out).

In the US, I acquired most of my mathematical education until I attended university classes from my parents. There was a designated chair in which my father would sit by my desk to explain new interesting math to me, and I was officially home-schooled in math for 2 years by him: I might add that those two years probably marked the peak of my teenage rebellion, but he patiently stuck it out, and one lesson I have held on to since those times was to listen to someone explaining some mathematics to me no matter if I thought I already knew it: chances are, my father said, I might learn something new. He was right then and is still right now: I~never dismiss mathematics as uninteresting and always try to learn something from it.

One might think that our close mathematical relationship through my childhood would mean a continued close mathematical relationship in adulthood. That has not been true: in college I~contemplated majoring in a scientific field besides math and leaning into it, out of a worry that mathematicians might think any of my success was due to my father's support. I remember talking about things along these lines to my father then, and he would often tell me, ``your problem is not your mathematical ability. It's your insecurity in it.'' It took me a few semesters in university to see that he was right, and I found a passion in number theory, which was fairly far removed from his mathematics (although he has told me at least 50 times that, while he would like to say he has never worked in number theory, that would not be entirely true). Since then our mathematical relationship has been much less intense, although I do still have a dream to produce a Fuchs--Fuchs paper sometime.

I have to add that, part of my father's dedication to my education in math certainly came out of his genuine care for teaching: even when he was teaching students in the university and not his children, he would go above and beyond the norm, writing his own lecture notes for nearly every class, shaping his own point of view in a way that I think is rare among university instructors. But part of it is a reflection of how much he cares for his family. His passion, I think, is split between two things: family and mathematics. For every fond memory that I have of us doing math together, I have one where math was not at all a factor. Putting up a tent together at countless campsites, hiking dozens of trails together, confiding in him both my biggest joys and fears. Most recently, I have seen the joy he gets out of seeing his grandchildren -- my two kids -- whenever they are around. I feel truly lucky to have him as my father, and wish him the most wonderful of birthdays, and many more to come.

\section{Ekaterina Fuchs}

\AuthorNameForHeading{Ek.~Fuchs}

It seems that the Russian way of raising a baby when I happened to be one involved a lot of walking, baby in stroller, outdoors. And so we walked, my Dad and I for hours on end. I would do what babies do, I suppose, and sleep, while he pushed the stroller, and did what he continues to do constantly~-- think and read. He would tell me how when it rained he would put a book into a clear plastic bag so he could read it without the book getting wet, and we would walk on, day in and day out.

And then night would fall, and he would read again, or rather recite, poetry to me until I~was soundly asleep. My mom has often told me stories of sitting in the kitchen (how vividly I~remember the layout of our 3~room apartment on {\it Trinadtsataya Parkovaya}, ``Thirteenth Park St'' in Moscow!) and knowing exactly how long it would be until he would be finished, since he recited the same poems every time, in his slow, peaceful baritone. Even when my sister and I~were older, he would often sit with us as we drifted off into sleep, reciting poetry to us.

\looseness=-1 Poetry continues to be a key motif in my memories of growing up and my dad. One summer, if I had to I would guess the summer of 1988, but I can't be terribly sure, my father and I decided together that he would help me memorize the entire anthology of Boris Pasternak poetry that followed the novel ``Doctor Zhivago''. Some days we would work on one poem, sometimes more than one, some poems stuck almost like honey in my mind, some were very difficult to internalize. The poems dealt with some very adult themes, and I remember in many cases our work to memorize the poems had to include long conversations about what the poem meant, and what the imagery was trying to evoke. We made it through that entire anthology that summer, just like we planned. To this day we occasionally allow ourselves a stroll down memory lane, and see if we can dust off this that or the other poem - while they are rusty and sleepy in my head, they are unbelievably crisp and alive in his. So many times I've said ``I'm going to relearn this one!'' and he's immediately ready to help me with any word or line where I stumble.

\begin{figure}[ht]\centering
\includegraphics[width=4.9in]{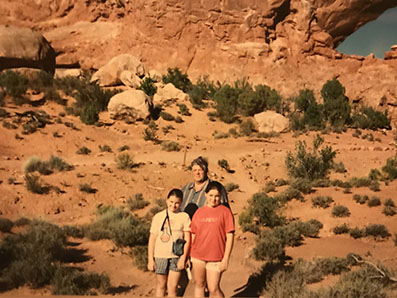}
\caption{On a hike with Katia and Lyalya.}
\end{figure}

I suppose this reminds me of another motif, if one could call it that, in my life with dad. The time he was willing to devote to basically anything we did together. One year he decided that the whole family needed a good solid lesson on the history of the Roman Empire. He immersed himself in historical texts and made copious notes, which he would then use to teach us (this included myself, my sister, my mom, and the cat). These were incredibly cozy evenings, sitting on the big L shaped sectional in the dining room, the three of us (and the cat) huddled together on the long side of the L, him sitting on the short end. I must confess that the evenings got so cozy, that my sister and the cat and I would often doze peacefully, which didn't go over particularly well with our parents.

He was then and continues to be now, unquenchably thirsty to teach us things -- whether it is beautiful mathematics, or the history of Ancient Rome, or poetry.

As a teacher now, I know that my desire to help others learn about beautiful things came from the true great fortune of having been surrounded by incredible teachers my whole life. My first teacher, without a doubt, was my dad.

I can honestly say that my dad's influence can be found at the very most foundational core of who I am today. His sense of humour, his dedication to his family, his love of art and thirst for knowledge, have shaped me in incredibly important ways. Basically, he's awesome, and I~wish him a happy birthday, and many, many, many more still to come.

\section{Alexei Sossinsky}

\AuthorNameForHeading{A.~Sossinsky}

I have known Dmitry Fuchs since 1957, when we were both students at the mathematics department of Moscow State University. We became neighbors quite by accident (our parents had acquired apartments in adjacent houses in Izmailovo) in our graduate student days and soon became close friends. This friendship has remained, despite the fact that, since 1990, we are usually separated by an ocean and a couple of continents.

Mitya, as everyone calls him is, above all, a mathematician, and while mathematics dominates his life, he has many other interests and activities that shape his unique personality: his sense of humor, his fascination for the outdoors, his love of poetry, his interest in the popularization and teaching
of mathematics (a rare trait for a research mathematician of his caliber). It
is mostly on these not very mathematical aspects of Dmitry Fuchs' life that I~intend to concentrate here, in the form of ``little Fuchs stories'' drawn from my vast oral repertoire.

\begin{figure}[ht]\centering
\includegraphics[width=3in]{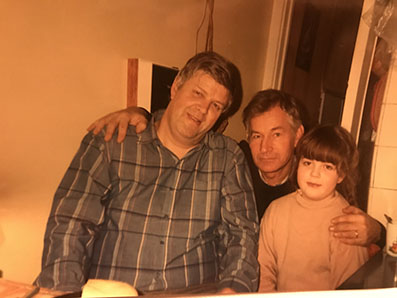}
\caption{D.~Fuchs, A.~Sossinsky, and Katia Fuchs.}
\end{figure}

{\it Sense of humour.} Mitya's is always low key, often ironic, it permeates his mathematical lectures and his conversation about people. Here is a related little story. We~-- a group mainly consisting of mathematicians~-- were downhill
skying on the Chiget mountain in the Caucasus. Yulia Gippenreiter, the non-mathematician in our group, was a psychologist and an accomplished (professional-level) Alpine skier. She would generously give us tips on skiing techniques; to Mitya, speeding fearlessly downhill, legs spread far apart for better balance, she would shout: ``Mitya, parallel skis, keep your skis parallel!'' Later, over tea in the cozy Ai chalet, Mitya a~would grumble, ``These stupid people in the humanities~-- they think that if your skis are more than a meter apart, they can't possibly be parallel'', and we would all burst out laughing.

\looseness=-1 {\it Poetry}. We first discussed poetry with Mitya as neighbors, in Izmailovo, and I must admit that at first~-- with a youthful snobishness of which I am now ashamed~-- I looked down on his knowledge and taste in Russian poetry. The first ``poetry reading'' where I began to change my mind occurred in the Kolsky peninsula in the late sixties, and I have described it in detail elsewhere, \cite{Tab},
but it still often comes to mind: the three of us (the third was our friend and colleague Lenya Erdman) on the Lovozero plateau, standing in a snowstorm near our tent, which was \mbox{being} progressively buried in snow, waiting for the blizzard to abate and the late Arctic morning to come, and~-- yes, you have guessed correctly~-- reading poetry to each other, sometimes in two or even three voices. I remember that Mitya got things started by reading a beautiful poem about Ulysses by Lugovskoi (whom until then I had regarded as a worthless Soviet establishment poet), Lenya recited a lot of funny poems by Sapgir and Zakhoder, I remember reciting Lermontov's {\it Demon} and, when my Russian repertoire had dried up, Poe's {\it Raven}, with my friends pitching in, like a Greek chorus, shouting out the mandatory ``Nevermore'' at the end of each stanza.

A less dramatic, but just as rewarding, reading took place a couple of years ago at the IHES in Bures-sur-Yvette, where we lived in neighboring cottages (not by accident, but by design). My present for Mitya's 80-th birthday was a selection of poems by Georgy Ivanov, a Russian \'emigr\'e poet he was not familiar with, and he responded with a short letter with a very precise literary analysis of Ivanov's verse, as good as those ever written by ``people in the humanities''.

\begin{figure}[ht]
\centering
\includegraphics[width=2.7in]{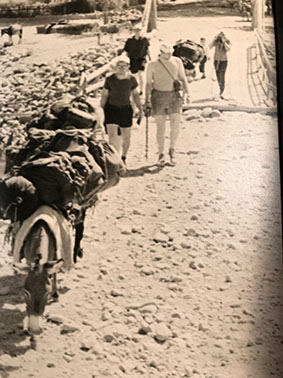}
\caption{Fuchs and Sossinsky follow the donkeys.}
\end{figure}

\indent{\it The outdoors.} Mitya could satisfy his love of nature only within the limits of the ``iron curtain'' surrounding the USSR, camping along its rivers, in its forests and mountains. This he would do within the framework of what we call {\it sportivny tourism}: whitewater trips (shooting rapids in kayaks), mountain hiking (on foot in summer or on skis in winter), hiking-camping trips in the rich forests of Russian midland. Of the several camping trips that we took together, I will describe the only one during which I took photographs: a high mountain hiking-camping trip through the Fanskie mountains, from Dushanbe to Samarkand. It began from the northern outskirts of Dushanbe, from which we began the long climb up with our backpacks loaded on two donkeys (see photo), until the path became so steep that the poor overloaded donkeys refused to continue. Progressively getting used to the increasing altitude, we hauled ourselves and our backpacks up to the first mountain pass, Mura, at 3500 meters.

\begin{figure}[ht]\centering
\includegraphics[width=3.2in]{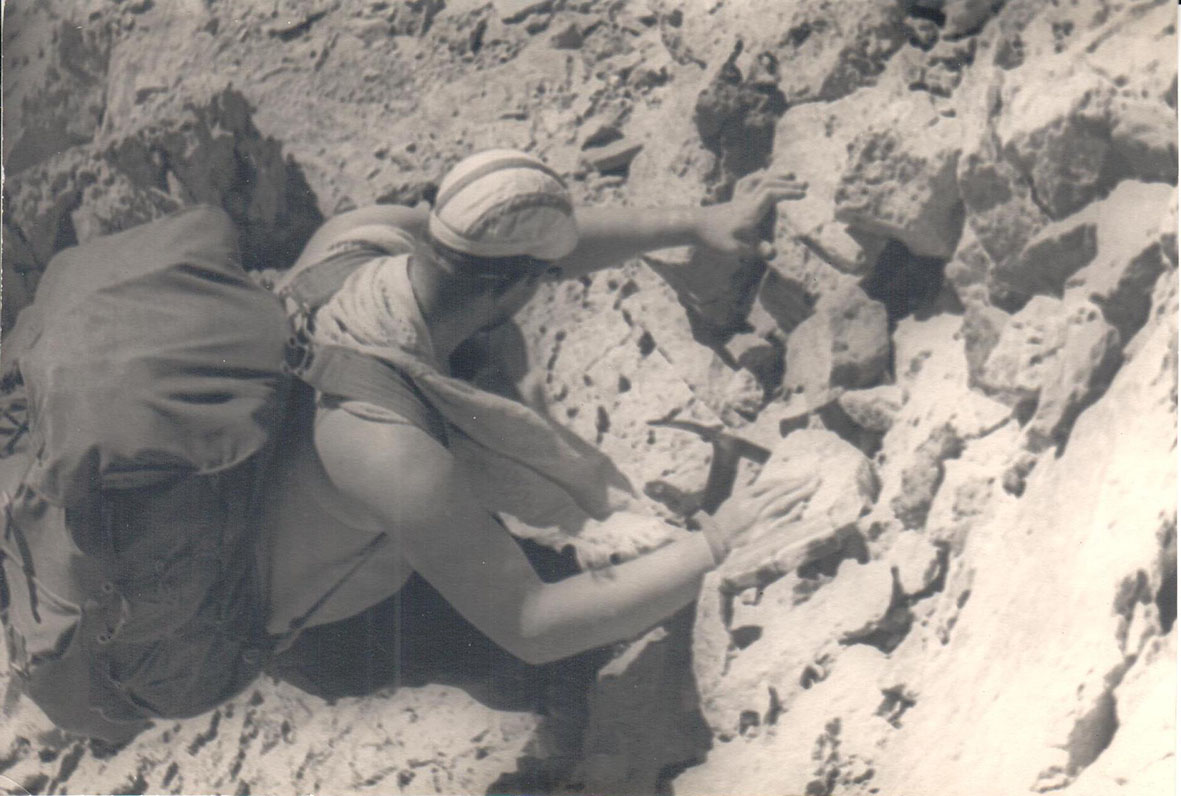}
\caption{Fuchs rock climbing.}
\end{figure}

From then on, it was down and then up over higher and higher mountain passes (the highest was at 4900 meters), crossing ice-cold streams, overcoming steep glaziers and rock climbing. Another photograph shows the group resting at the start of one of the glaciers, with two local shepherds gratefully sharing our meal. The local population (nomadic shepherds) was very sparse and practically untouched by 20th century civilization, their self-sufficient way of life unchanged for many centuries, most of them having never handled money, been in a city, or seen a photo-camera; they looked at us as if we were creatures from another planet.

\begin{figure}[ht]\centering
\includegraphics[width=3in]{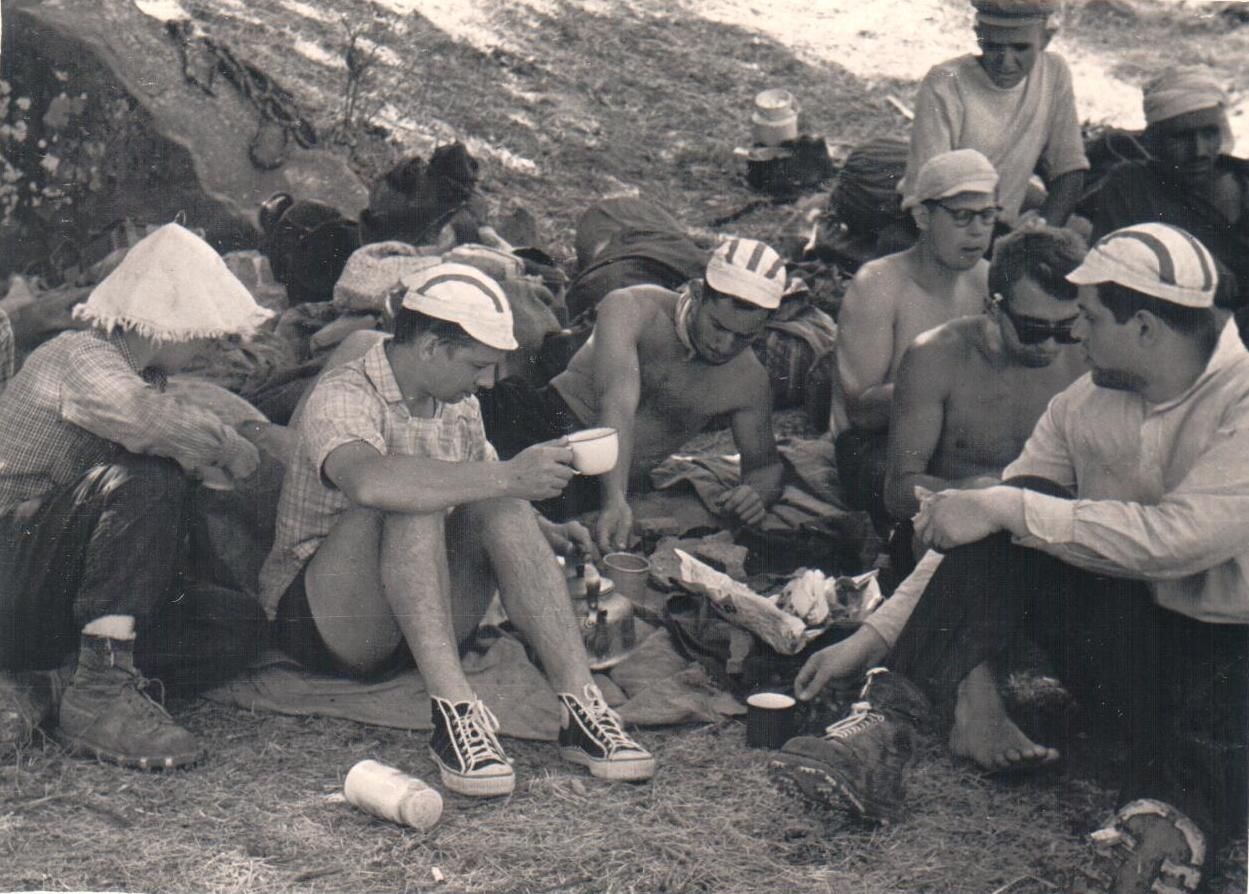}
\caption{A frugal meal at the foot of a glacier.}
\end{figure}

One of the high moments of our trip (literally and figuratively) was in the morning at our camp site some 4000 meters above sea level, when a huge eagle swooped down on us, grabbed a~sausage that Mitya was about to cut up into seven equal pieces and flew off with it; I think it was then that Mitya philosophically declared: ``Well, the sausage was beginning to spoil anyway''.
The whole episode took two or three seconds, and so, unfortunately, there is no photograph of the eagle, and there are no photographs that do justice to the mountain scenery and the beautiful and then almost totally uninhabited Iskander Kul' lake (now a tourist resort), on whose shores we rested for a couple of days before finishing our trip in the historic city of Samarkand.

\begin{figure}[ht]
\centering
\includegraphics[width=10cm]{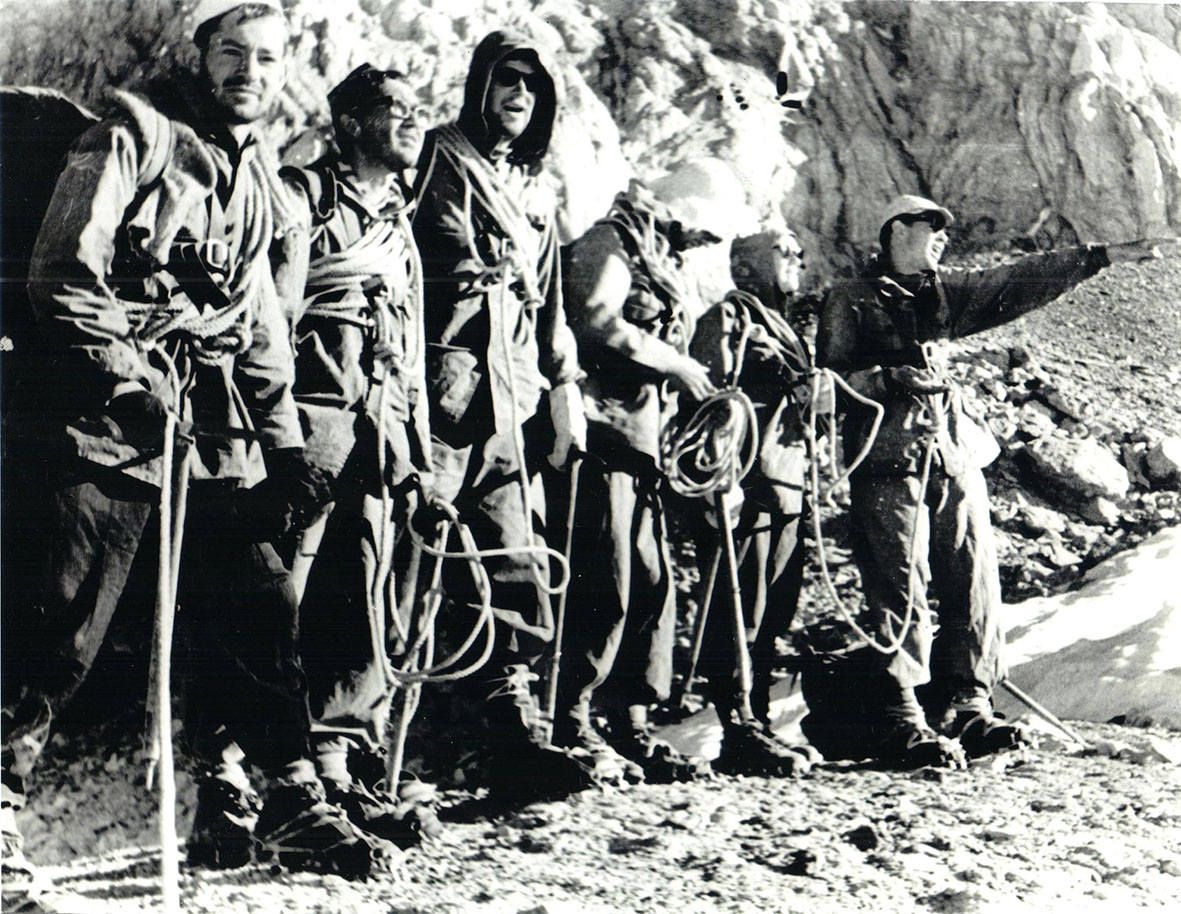}
\caption{ A.~Sossinsky, I.~Girsanov, D.~Fuchs, E.~Mochigin, G.~Turina, Yu.~Chapovsky ({\it left to right}).}
\end{figure}

There is, however, one photograph that is really worth looking at. It very expressively shows six tired climbers, some still breathing heavily, who have just overcome a mountain pass and look forward with optimism to the continuing trip. Yet this is a tragic photo: five years after it was taken, only two of the six persons on the picture were still in this world: Girsanov, Mochigin, and Chapovsky were buried in an avalanche in the Sayany mountains in 1967 and Galya Turina, Mitya's wife, drowned two summers later when her kayak overturned in a white water trip in 1970 just above the Arctic Circle in the Urals. But Mitya and I are still here and look back on this trip as one of the most rewarding experiences of our lives.
\medskip

\centerline{* * *}

\medskip

Perhaps I should conclude this text by shortly recalling how our life lines converged and diverged in the last 60 years. From our graduate years on, our lives first developed in parallel, with Mitya one year ahead of me (although almost two years younger, he was a year ahead in his studies). After completing our graduate studies and defending our PhD's, we became, one after the other, associate professors at the topology Chair of P.S.~Alexandrov (although neither of us were his pupils, and we worked in the rapidly developing field of algebraic topology that Alexandrov no longer understood). We both became regular participants in the famous Gelfand seminar, and also worked together in the math olympiad movement on the national and international level.

I recall an episode that occurred at the Gelfand seminar: during a break, I walked up to the blackboard and tried to explain to Fuchs how to answer a question raised at the talk preceding the break; this involved the calculation of the first term of a spectral sequence, I explained the filtration and, with Mitya's help, indicated how the calculation would go. I did not notice that Gefand was standing by and attentively listening; when I had finished and the seminar was about to resume, Izrael Moiseevich, addressing us both, said that was quite interesting and asked us to write it up and bring the text next week, he would have it published as a joint paper by the three of us. But we didn't: Mitya convinced me that the question was not very serious and the proof too simple to be worth publishing.

Perhaps a year later, at the end of the seminar, Gelfand took us aside and proposed us research positions at his mathematical biology laboratory. In my case, I think it was the episode described above which decided that I was worth hiring; as to Mitya, it was by then clear that he was a first rate mathematician who would embellish any research center. Mitya accepted. I~refused: I~was simply afraid of Gelfand.

This was a second step in our moving apart, the first one being geographical:
I moved to an apartment in the southern part of town, so we were no longer neighbors. The next step in that direction occurred when I was forced out of Moscow University, and after a year of unemployment, started working at the popular science magazine {\it Kvant}. During the 13-year period that I worked there, we would meet with Mitya at the {\it Kvant} headquarters more often than at research seminars. Mitya became a regular contributor to the magazine, and an unexpected result of that was his marriage to Ira Klumova, who then worked as a {\it Kvant} math editor. It was my wife, Elena Efimova, not I, who played the role of the matchmaker: two side-by-side tickets to a theater performance given by Lena to Mitya and Ira started things off \dots

Another place where we continued to meet was Bella Subotovskaya's {\it People's University} (while it lasted), where Mitya, in this case following my footsteps, became an active teacher until this ``university'' was closed down by the KGB. I have described the dramatic story of that unusual institution in my article~\cite{Sos}.

One of our last contacts in Moscow before Mitya's emigration to the US in
1990 was a series of talks on knot polynomials that I gave at the Fuchs--Varchenko seminar. This was extremely important for me~-- Mitya's very favorable reaction and encouragement helped me realize that during my ``exile'' to {\it Kvant} I did not lose the ability of doing serious mathematics.

\looseness=-1 Since 1990, our meetings have been few and far between: at Davis, of course,
at Oberwolfach, at the IHES, in Moscow at the Gelfand Centenary celebration~-- are those that first come to mind.

Birthday wishes? I have sent them to Mitya by e-mail: to be able to think clearly and mathematically till the very end, to derive from life whatever pleasures we are still capable of appreciating, and, our healths and God willing~-- a~few more meetings on either side of the ocean.

\section{About Mitya Fuchs. By Alexandre Kirillov}

\AuthorNameForHeading{A.~Kirillov}

I was asked to write something for the special issue in honor of 80th anniversary of D.B.~Fuchs, my old friend, an outstanding mathematician, and simply a very good man. Of course, I~immediately agreed. But when push came to shove, it turned out that to do it wasn't that easy. In general, it is not easy to write a good mathematical paper, and to write it better (or at least not worse) than one's previous works is even harder. Let alone on a short notice \dots

Of course, one could simply write about Mitya, how we lived, worked, what we did in our free time. But this immediately forces one to remember many other people and all our life in the USSR, Russia, and America because in all of this Mytya played an important role. And still, I shall give it a try.

\begin{figure}[ht]\centering
\includegraphics[height=2.3in]{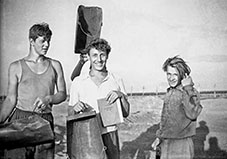}\ \ \
\includegraphics[height=2.3in]{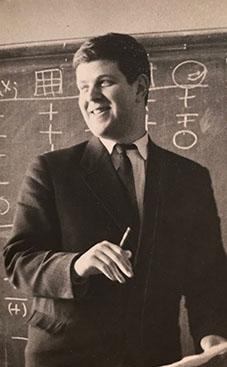}
\caption{Left: as part of Khrushchev's campaign to open vast tracts of virgin land in the northern Kazakhstan and the Altai region of the USSR, the students were sent to work there in summer. Right: a young topologist, D.~Fuchs.}
\end{figure}

Much was said about the Golden Years of Moscow mathematics (1955--1967). We both had the privilege
to study at MekhMat of MGU (Mathematics Department of the Moscow State University) that, without exaggeration, was then the leading mathematical center in the world.

\looseness=-1 Among our professors were P.S.~Alexandrov, N.V.~Efimov, I.M.~Gelfand, A.O.~Gelfond,
B.N.~Delone, E.B.~Dynkin, A.N.~Kolmogorov,
A.G.~Kurosh, A.A.~Markov, D.E.~Menshov,
I.G.~Petrovski, L.S.~Pontryagin, M.M.~Postnilov, I.R.~Shafarevich, B.V.~Shabat, G.E.~Shilov. The
recitations and research seminars for youngsters were run by then young, but already well known,
mathematicians F.A.~Berezin, R.L.~Dobrusin, R.A.~Minlos, A.S.~Schwarz, A.G.~Vitushkin.

The general atmosphere of those years (1955--1967) was
splendid, and it has never been the same again. What
happens at MekhMat today does not stand comparison.
I was lucky to belong to a group of 20--25 students who
were united in their admiration for talent and were
genuinely interested in each other's work. (To mention
only the most active and well remembered:
A.~Arkhangelskij, V.~Arnold, B.~Averbuch,
N.~Brushlinskaya, A.~Chernavskij, L.~Churaeva, V.~Murskij,
R.~Sofronitskaya, M.~Shur, V.~Tutubalin, F.~Vetukhnovskij,
E.~Vinberg.)

Later our group extended by students from the next year
(D.~Fuchs, V.~Palamodov, A.~Sossinsky, G.~Tyrina) and
some previous years (V.~Alexeev, I.~Girsanov, B.~Polyak,
Ya.~Sinai, V.~Tikho\-mi\-rov).

Besides mathematics, the common trend was
tourism. Not in the western form of expensive tours in
comfortable conditions, but week-end pedestrian hiking
on foot near Moscow, or longer travels during student
vacations on ski, on kayaks, in the mountains. These
travels were partially sponsored by the University Tourist
Club. One can get some outfit (skis, kayaks, ropes) and
some foodstuffs (condensed milk, braised meat) which
was rather difficult to find in the post-war Russia.

Mitya was a very important member of our expeditions
due to his strength, endurance, and good nature. Once he
saved our kayak when we tried to pass in poor visibility through
a narrow passage in a dam. The kayak was almost folded in
half, being pressed by the waterfall and only Mitya's
power legs and spine allow us to pull it out.

Another remarkable situation happened during our
winter expedition on the Northern Ural, the coldest place
in European Russia, where the temperature in winter falls
to $-50$~C$^\circ$. By the way, one of our common discovery with Mitya
 was that, despite the general opinion, the spittle
(and other warm excretions) do not freeze immediately in
the air even under this temperature.

\looseness=-1 Preparing to this expedition, we decided to make an
experiment and take with us a portable stove working on
firewood (the gas stoves did not yet exist in tourist
practice at that time). On the third or fourth day, when the
frost became rather severe, we installed the stove in our tent
and started the fire. Soon the air inside the tent became
rather warm, but the whole space filled up by the dense
smoke. Since none of us had experience of heating a tent
with a stove, we thought that it was a natural effect and
tried to endure the smoke. But after a short time the smoke
turned into flame, and when we finally succeeded to put it
out, almost half of our tent had gone.

The culprit of the fire was the veneer shovel,
incautiously leaned against the stove. (Such shovels are
usually used in winter expedition for digging the snow to make
a place for the tent and for the campfire.) In the morning
there was a council about one main question: what to do?
Some people proposed to fix the tent and to go further. Other
preferred to break the expedition and come back home.

Since I was the official head of expedition, my vote was
decisive. And with great doubt I~decided to return. I was
very thankful to Mitya and his wife Galya who did not
blame me for this decision, though I suspect they felt
differently. To my surprise, the way back was rather easy
and agreeable. First, going down the river and using our
trace in snow, we made the way back in one day (more
precisely, in 12~hours). Second, the frost was softened and
when we came to the railroad station, the temperature was
$-18$~C$^\circ$, feeling like a summer day.

The mathematical achievements of Mitya are numerous
and well-known. Last several years we tried to work
together. Luckily, I draw to this project my former PhD
student V.~Ovsienko and his wife S.~Morier-Genoud. Thanks to Valya's administrative talent and
modern social trends in scientific politics, the four of us formed
a ``Square'' team at the American Institute of
Mathematics (Palo Alto~-- San Jose).

We spend at AIM three
one-week terms, and wrote the joint article, where we
tried to realize my idea about the application of the orbit
method to the group of triangular matrices over a~finite
field. I'll explain here the problem in very informal way.

\looseness=1 The point is that usually the linearized version of a
problem is simpler than the non-linear original. But the
classification problem for coadjoint orbits of the
triangular matrix group $N$ is still open, notwithstanding
that it is the linearized version of the classification
problem for the Borel group orbits in the flag manifold $F
$. The solution of the latter problem is well-known: the
B-orbits are labeled by the elements of the Weyl group. I
conjectured that the coadjoint orbits are in some sense the
``shadows'' of $B$-orbits. The naive formalization of this
conjecture turns out to be wrong, but we have found a~very interested new phenomena and are eager to continue.

I can not speak about Mitya and not to mention his family.
His first wife, Galina Tyurina, was an ardent tourist and
a downhill skier. She was also an outstanding mathematician and
an essential member of I.R.~Shafarevich's team in
algebraic geometry. She died tragically during a~kayaking
expedition.

The second Mitya's wife, Irina Klumova, was my student
at MechMat and later my co-author in the magazine
``Kvant'', where she worked for many years. They have two
beautiful daughters, Katia and Elena, both talented
mathematicians. Thanks to Elena, Mitya became a~happy
grandfather of two.

\begin{figure}[ht]\centering
\includegraphics[height=2.3in]{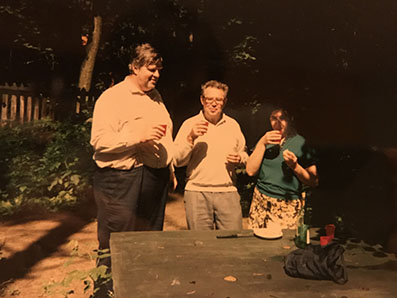}\ \
\includegraphics[height=2.3in]{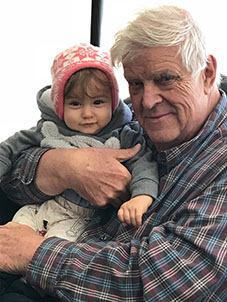}
\caption{D.~Fuchs, A.~Kirillov, and I.~Klumova. Fuchs, the grandfather.}
\end{figure}

Now we meet not as often as we wish, but I still look forward in the hope on the further joint work.

\section[Close friendship with the Fuchs family in Moscow. By Alice Fialowski]{Close friendship with the Fuchs family in Moscow.\\ By Alice Fialowski}

\AuthorNameForHeading{A.~Fialowski}

When I was a graduate student (an aspirant) in Moscow from 1979 to 1983, I lived in the dormitory of the Moscow State University, MGU. I regularly visited the Gelfand Seminar on Monday evenings. In a sense, it was a big family in the lecture room: students and professors were talking to each other, made friends there. One of the main people at those seminars was Dmitry Fuchs, who usually sat with Alexandre Kirillov in the first row. Gelfand often pestered them with questions.

My first year as an aspirant was quite hard. Partially it was related to the fact that I did not have much background in the topics which were discussed. MGU was the very best school of mathematics in the entire world at the time, and the research topics there were numerous and at the highest international level. In comparison, the Hungarian school was mainly influenced by Erd\"os and combinatorics in those days. I liked algebra and functional analysis, and started studying Lie groups, Lie algebras, and representation theory. I was also learning Lie algebra cohomology from Fuchs, and from his famous book which was already available in Russian. In order to use this theory for my project in algebras, I needed help.

\begin{figure}[ht]
\centering
\includegraphics[height=2.2in]{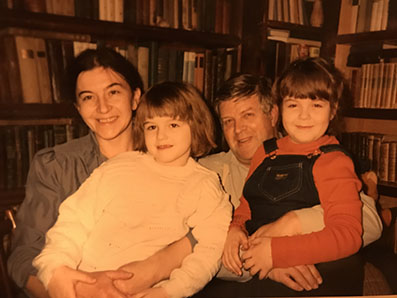}\quad
\includegraphics[height=2.2in]{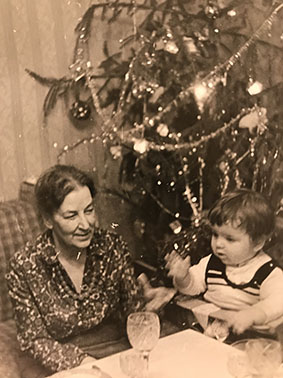}
\caption{The Fuchs family. Ekaterina Ivanovna.}
\end{figure}

In those days in Moscow professors did not have offices and used to work with their colleagues at home. Fuchs invited me to his home and I went to visit him and his mother, Ekaterina Ivanovna. We had dinner together, wonderful Russian salads and other stuff. I was touched by their openness and hospitality. We worked for a few hours, decided to meet the next week, then the next, and thus our Wednesday meetings became regular.

Of course, we discussed many subjects: the Fuchses were interested in my country, my family, and Dmitry Borisovich explained to me some details about the situation in Russia, talked about his brother, his father, his family. I met Ira, his wife, I was there when they had their first baby Katia, and also when their second baby, Lyalya, was born. I also went with them to their dacha. I will forever be grateful to them for making me feel a part of their family, and the Wednesday meetings were really important to me. Living alone in Moscow in those days as a foreigner was not so easy, but this way I got a family there!

I became especially close with Ekaterina Ivanovna, Fuchs's mother. She was a lovely lady, and there was something warm and reassuringly peaceful about her. When I was depressed (which unfortunately happened from time to time), she was always able to help me get away from feeling down. I was also there when she was dying. She died in her family, in peace.

After my thesis defense, I returned home to Budapest, but I remained in regular contact with my ``Russian family''. Fuchs, Ira, and the girls visited me in Budapest, and I was able to take them to nice playgrounds on the Gell\'ert-hegy, and other places. Dmitry Borisovich had many good friends in Budapest, like
K\'aroly M\'alyusz, Andr\'as Kr\'amli, P\'eter Major, Doma Sz\'asz, and others. They all got acquainted with Dmitry Borisovich in Moscow, and his lovely personality attracted many-many people. (Unfortunately K\'aroly and Andr\'as are not with us anymore \dots)

\begin{figure}[ht]\centering
\includegraphics[width=2.6in]{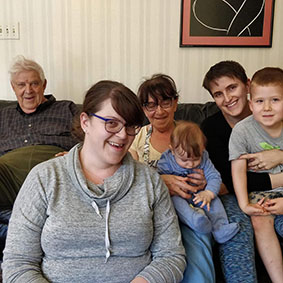}
\caption{The Fuchs family.}
\end{figure}

Soon after I returned home, I began to get invitations to the West. I was able to go to Germany for a Humboldt fellowship and to attend various conferences. I also received an invitation to the University of Pennsylvania and later I got a permanent job at the University of California at Davis. While I was there I learned that Dmitry Borisovich might be able to come to the United States with his family. I immediately suggested to our Department Chair to invite him to join our Faculty. Eventually, it worked out and since then Dmitry Borisovich Fuchs, together with Albert Schwarz, have graced the Davis mathematics community. UCDavis is lucky to have them! Recently Fuchs's younger daughter, Elena also joined the Faculty, together with her husband, Martin Thanh Luu.

I was able to continue working with Fuchs in Davis as well. We lived just a few houses apart, and I often visited the Fuchs family there too. After my return to Hungary, we have not seen each other very often. But I went to Davis for his 70th birthday celebration, and we also met every other summer at the Max-Planck Institute in Bonn, where both his family and mine usually spend one month.

I feel lucky that I have such a real friend as Dmitry Borisovich. For me he is not just a~colleague with whom it is good to work, but a close friend, with whom we share a~lot of values in life. I wish him all the very best and long and happy years with his family!

\section{Personal notes by Valentin Ovsienko}

\AuthorNameForHeading{V.~Ovsienko}

What can be added to the beautiful (and sometimes romantic) descriptions of the mathematical life in Moscow
just before the {\it Perestroika}?
Many remarkable memories have been published in recent decades by various authors.
``Those days are past now and in the past they must remain~\dots''
I can only talk about personal impressions of a young undergraduate student that I was at that time.

The real life at {\it Mekhmat} started in the afternoon, and I spent them all
 running from one fascinating seminar to another.
This is how I met Dmitry Borisovich.
I saw him several times every week, at Gelfand's seminar
(that, as we know, was the main mathematical event), Arnold's seminar, meetings of the Moscow Mathematical Society, etc. Fuchs gave several lectures at Math Society,\footnote{A series of survey introductory lectures for a large audience.}
where he talked on various subjects.
I remember very well his lecture about symplectic topology, a new domain
that was making its ``first steps''.
It was not at all (not yet) among the subjects of his own research,
but Fuchs gave a very clear and precise overview.
I was not at all surprised, it was obvious to me that Fuchs knew everything!

My first personal meeting with Dmitry Borisovich was related to the seminar of Arnold.
I~was a regular member of the seminar, and Arnold invited me
to give a talk about the ``Lagrange Schwarzian derivative'', a notion introduced in one of
my first research papers.
That was my first talk outside Kirillov's seminar for younger students,
and I was very anxious.
Being out of town, Arnold asked Fuchs to run the seminar, and Dmitry Borisovich took the job seriously:
he suggested to meet and discuss a few days before the seminar.
Our meeting became for me one of the episodes that remain in memory for the rest of the life.
Fuchs asked many questions and made several comments, that later I used many times.
Trying to explain my work to him, I got an impression of understanding it myself.
My talk at the seminar went smoothly \dots

\looseness=-1 My research was largely influenced by the famous papers authored by Gelfand--Fuchs and Feigin--Fuchs.
With varied success, I tried to connect their cohomology classes to projective differential geometry,
and this task occupied me for many years.
During these long years I met Dmitry Borisovich from time to time at mathematical conferences and, even more often, at the Tabachnikovs' home.
But never could I imagine that I would collaborate with him,
until such an opportunity suddenly presented itself.
This happened 30 years after the events described above.

\begin{figure}[ht]\centering
\includegraphics[width=2.0in]{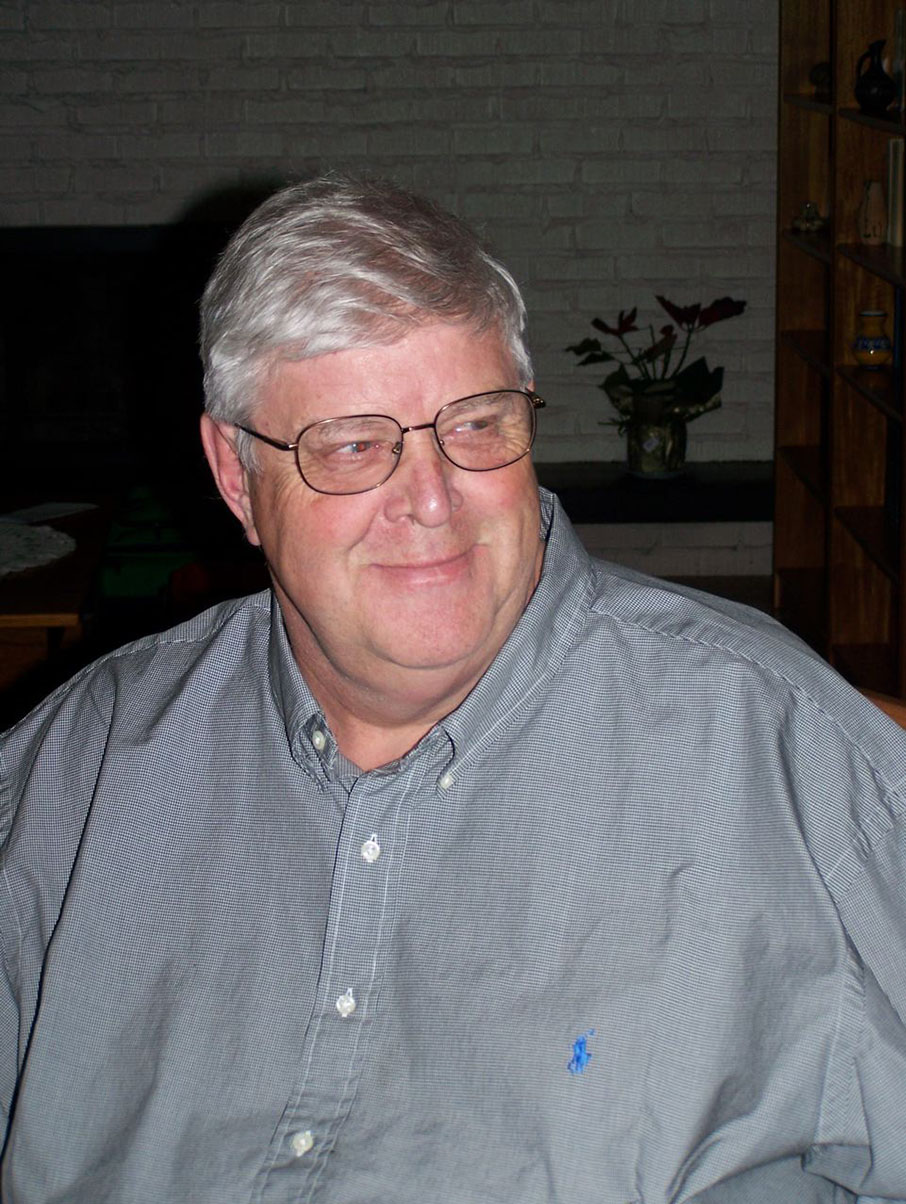}
\caption{Dmitry Borisovich Fuchs.}
\end{figure}

Five years ago, my teacher, Alexandre Alexandrovich Kirillov, who has always been a close friend of Dmitry Borisovich, suggested to work on an old unsolved problem of classification of coadjoint orbits of the nilpotent Lie algebra of (upper) triangular matrices. Kirillov created a~``square-team'' consisting of four people: Fuchs, Kirillov, Sophie Morier-Genoud, and myself. We worked within the SQuaREs program at AIM, and later, at MFO, Oberwolfach. Let me mention that during our work on the project that lasted for three years and a half, the square eventually became a pentagon.

This work with my teachers was yet another moving experience. Fuchs and Kirillov stayed at the blackboard, intensively discussing the subject and paying almost no attention to us. Kirillov directed the general attack, Fuchs generated the technique. Sometimes they asked Sophie about details of the modern combinatorial theories. I was the chronicler, always with a notebook at hand, and I recorded each word of the ``oracles".

In the afternoon, the discussion usually became more relaxed.
``Do you remember when we climbed Tian Shan Mountain in 1956?'' would ask one of them.
``Of course, I remember'', would reply the other, ``but it was in 1957!''

Gradually, our subject evolved.
Approximating the orbits, we found ourselves on a more solid soil, studying Schubert varieties and their tangent cones.
Little by little, some understanding and some results started to appear.

However, writing of our paper \cite{MR3766071} was not an easy promenade.
The first version contained two theorems that would solve our problem completely.
I was in charge of reconstructing the details of the proofs, outlined on the blackboard by the collective efforts, and spent many sleepless nights on that.
Finally, I wrote an SOS message to my coauthors, asking for a permission to call the second theorem a ``conjecture''.
When permission was granted, Sophie found a counterexample to this brand new conjecture.
And, finally, the ``only if'' part of the first theorem also became a~conjecture, but fortunately, with no counterexamples in sight.

\pdfbookmark[1]{References}{ref}

\end{document}